\documentclass[letterpaper,11 pt,conference]{ieeeconf}
\usepackage{graphicx}
\usepackage{amssymb}
\usepackage{amsmath}
\IEEEoverridecommandlockouts

\newtheorem{remark}{Remark}
\newtheorem{definition}{Definition}
\newtheorem{theorem}{Theorem}

\title{Certifying non-existence of undesired locally stable equilibria in formation shape control problems}
\author{Tyler H. Summers, Changbin Yu, Soura Dasgupta, Brian D.O. Anderson
\thanks{T.H. Summers is with the Automatic Control Laboratory, ETH Zurich, Switzerland, \{tsummers@control.ee.ethz.ch\}. C. Yu and B.D.O. Anderson are with the Australian National University and NICTA, Canberra, Australia, \{brad.yu,brian.anderson\}@anu.edu.au. S. Dasgupta is with the Department of Electrical and Computer Engineering, University of Iowa, USA, \{dasgupta@engineering.uiowa.edu \}.%
 }%
 \thanks{This work is in part  supported by the ETH Zurich Fellowship Program and by the US NSF grants CCF-0830747, EPS-1101284  
and CNS-1329657, ONR grant N00014-13-1-0202 and a grant from the Roy J. 
Carver Charitable Trust.}%
}

\begin{document}

\maketitle

\begin{abstract}
A fundamental control problem for autonomous vehicle formations is formation shape control, in which the agents must maintain a prescribed formation shape using only information measured or communicated from neighboring agents. While a large and growing literature has recently emerged on distance-based formation shape control, global stability properties remain a significant open problem. Even in four-agent formations, the basic question of whether or not there can exist locally stable incorrect equilibrium shapes remains open. This paper shows how this question can be answered for any size formation in principle using semidefinite programming techniques for semialgebraic problems, involving solutions sets of polynomial equations, inequations, and inequalities. 
\end{abstract}

\section{Introduction}
Autonomous vehicle formations have received significant attention in the systems and control literature recently, including theoretical research and applications in teams of unmanned aircraft for reconnaissance and surveillance missions, satellite formations for deep-space imaging, and underwater vehicle swarms for oceanic exploration. One fundamental control problem is formation shape control, in which the agents must maintain a prescribed formation shape using only information measured or communicated from neighboring agents. 

A substantial line of research in the literature involves distance-based formation shape control, in which agents move to maintain relative distances between certain neighboring agents \cite{krick2009stabilisation,yu2009control,summers2011control}. Due to the focus on distances, the research features tools from nonlinear systems and control theory and has strong connections to the theory of rigid graphs \cite{anderson2008rigid}. In \cite{krick2009stabilisation}, gradient-descent control laws were proposed and shown to provide local exponential stability of the desired formation shape whenever the graph of the underlying information architecture is infinitesimally rigid. 

Global stability properties for distance-based formation control problems remain a significant open problem. Some global stability results are known for three-agent formations \cite{cao2007controlling,anderson2007control,cao2011maintaining}; for example, in \cite{cao2007controlling} it is shown that for almost all initial conditions (except for a thin set of collinear initial conditions), the formation converges exponentially fast to the desired shape. For four-agent and larger formations, much less is known. It was observed in \cite{krick2009stabilisation} in a specific example in which all inter-agent distances are actively controlled that under their gradient-based control law, there can exist non-degenerate incorrect equilibrium shapes. For two-dimensional rigid formations of point agents, this turns out to be an inherent consequence of using any gradient-based control law \cite{anderson2011morse}. However, it was recently shown that when the desired formation is rectangular, the associated incorrect equilibrium shapes are locally unstable \cite{anderson2010controlling,dasgupta2011controlling}. 

Even in four-agent formations, the basic question of whether or not there can exist locally stable incorrect equilibria remains open. The present paper shows how this question can be answered in principle for any size formation using semidefinite programming techniques for semialgebraic problems, involving solutions sets of polynomial equations, inequations and inequalities. The techniques are based on the recent observation by Parrilo \cite{parrilo2000structured,parrilo2003semidefinite} that determining if a given polynomial is a sum of squares can be done efficiently by solving a semidefinite program, which is a convex optimization problem. Combining this technique with a fundamental result from real algebraic geometry, called the \emph{Positivstellensatz}, one can computationally certify if a given semialgebraic set, described by finite sets of polynomial equations, inequations and inequalities, is empty. For formation shape control problems using the gradient-descent control laws proposed by \cite{krick2009stabilisation}, the conditions for equilibrium, shape incorrectness and local stability  can all be described by finite sets of polynomials. Thus, these techniques provide a way to certify non-existence of locally stable incorrect equilibria in principle.

The rest of the paper is organized as follows. Section II briefly reviews semidefinite programming techniques for sum of squares polynomials and the Positivstellensatz from real algebraic geometry. Section III describes how the formation shape control problem fits into this framework and shows how these techniques can be used to generate infeasibility certificates in this context. Section IV gives concluding remarks. 

\section{Semidefinite Programming and Real Algebraic Geometry}
This section gives a very brief introduction to sum of squares programming and real algebraic geometry. In particular, we recall how semidefinite programming can be used to choose coefficients in a polynomial so that it is a sum of squares, and how the Positivstellensatz can be used to certify that a semialgebraic set is empty. Further details can be found in \cite{parrilo2003semidefinite} and \cite{bochnak1998real}.

\subsection{Sum of Squares Polynomials and Semidefinite Programming}
We begin with some basic definitions.
\begin{definition}
A \emph{polynomial} $f$ in $x_1,...,x_n$ is a finite linear combination of monomials
\begin{equation}
f = \sum_\alpha c_\alpha x^\alpha = \sum_\alpha c_\alpha x_1^{\alpha_1}\cdots x_n^{\alpha_n}, \quad c_\alpha \in \mathbf{R},
\end{equation}
where the sum is over a finite number of $n$-tuples $\alpha = [\alpha_1,...,\alpha_n]$, $\alpha_i \in \mathbf{N}_0=\{0,1,2,...\}$. The set of all polynomials in $x\in \mathbf{R}^n$ with real coefficients is denoted by $\mathbf{R}[x_1,...,x_n]$.
\end{definition}

A polynomial is called sum of squares if it can be expressed as a finite sum of squared polynomials. The set of sum of squares polynomials is denoted by
\begin{equation}
\begin{aligned}
 \Sigma = \{ F(x) \in \mathbf{R}[x_1,...,x_n] \ | \  F(x) = \sum_i f_i^2(x), \\
  f_i(x) \in \mathbf{R}[x_1,...,x_n] \}
\end{aligned}
\end{equation}
and forms a \emph{cone} in $\mathbf{R}[x_1,...,x_n]$. 

Let $F(x) \in \mathbf{R}[x_1,...,x_n]$ be a polynomial of degree $2d$. 
 Let $z(x)$ be the vector of all monomials of degree less that or equal to $d$. The following result is from \cite{parrilo2003semidefinite}.
\begin{theorem}[\cite{parrilo2003semidefinite}]
The polynomial $F(x)$ is a sum of squares if and only if there exists a symmetric positive semidefinite matrix $Q$ such that $F(x) = z^T(x) Q z(x)$. 
\end{theorem}

Thus, by choosing a positive semidefinite matrix subject to affine constraints defined by matching coefficients of $F(x)$ with an expansion of a quadratic form in the monomial vector $z(x)$, one can ensure that $F(x)$ is a sum of squares. This is a semidefinite programming feasibility problem in primal form, which is a convex optimization problem and can be solved efficiently using various well-developed software packages. This method can be used if the coefficients of $F(x)$ must satisfy some additional affine constraints.

\emph{Example}. Consider the quartic form in two variables: $$F(x,y) = a_1 x^4 + a_2 x^3 y + a_3 x^2 y^2 + a_4 x y^3 + a_5 y^4 $$ and suppose we require affine constraints $Ga = h$ on the coefficients $a = [a_1,...,a_5]^T$ to be satisfied.  The vector of monomials can be taken as $z(x,y) = [x^2, xy, y^2]^T$ since $F$ is a form (i.e. the total degree of all monomials is the same, here four). Now set
\begin{equation}
\begin{aligned}
F(x,y) &= a_1 x^4 + a_2 x^3 y + a_3 x^2 y^2 + a_4 x y^3 + a_5 y^4 \\
           &= \left[\begin{array}{c}x^2 \\xy \\y^2\end{array}\right]^T \left[\begin{array}{ccc}q_{11} & q_{12} & q_{13} \\q_{12} & q_{22} & q_{23} \\q_{13} & q_{23} & q_{33}\end{array}\right] \left[\begin{array}{c}x^2 \\xy \\y^2\end{array}\right]\\
           &= z(x,y)^T Q z(x,y)\\
           &= q_{11}x^4 + 2q_{13}x^3y + (q_{33} + 2q_{12})x^2y^2  \\
           & \quad + 2q_{23}xy^3 + q_{22}y^4
\end{aligned}
\end{equation}
Matching coefficients, we obtain the following semidefinite program with decision variables $Q\in S_+^n$ ($S_+^n$ denotes the set of symmetric positive semidefinite matrices) and $a \in \mathbf{R}^5$
\begin{equation}
\begin{aligned}
& \text{find} && Q \succeq 0 \\
& \text{subject to} && q_{11} = a_1, \quad 2q_{13}=a_2, \quad q_{33}+2q_{12} =  a_3, \\
&			    &&  2q_{23}=a_4, \quad q_{22} = a_5 \\
&			    && Ga=h.
\end{aligned}
\end{equation}
If the problem is feasible, the polynomial coefficients $a$ and factorization matrix $Q$ are returned. The problem can be infeasible if the additional affine constraints $Ga=h$ make it impossible for $F$ to be a sum of squares. 

\subsection{The Positivstellensatz}
This section describes the Positivstellensatz, a fundamental result from real algebraic geometry that establishes a relationship between the set of real solutions to a set of polynomial equations, inequations, and inequalities (a semialgebraic set) and algebraic concepts of polynomial cones, ideals and monoids. Combining this result with the semidefinite programming techniques from the previous section gives a computational method for determining if a given semialgebraic set is empty. We begin again with a few basic definitions.
\begin{definition}
The set $I \subseteq \mathbf{R}[x_1,...,x_n]$ is an \emph{ideal} if it satisfies:
\begin{itemize}
\item $0 \in I$
\item $a,b \in I \Rightarrow a+b \in I$
\item $a \in I, \ b \in \mathbf{R}[x_1,...,x_n] \Rightarrow a b \in I$
\end{itemize}
\end{definition}

\begin{definition}
Given a finite set of polynomials $(f_i)_{i=1,...,s}$, define the set
$$\langle f_1,...,f_s \rangle = \left\{ \sum_{i=1}^s f_i g_i \ | \ g_i \in \mathbf{R}[x_1,..,x_n] \right\}.$$
This set is easily shown to be an ideal and is called the ideal \emph{generated} by the $f_i$.
\end{definition}

\begin{definition}
Given a finite set of polynomials $p_i \in \mathbf{R}[x_1,..,x_n]$, the set $M(p_i)$ of finite products of the $p_i$ (including the empty product and the identity product) is called the \emph{multiplicative monoid} generated by the $p_i$. Note that powers of the $p_i$ are also allowed.
\end{definition}

\begin{definition}
A subset $P \subseteq \mathbf{R}[x_1,..,x_n]$ is called a \emph{cone} if it satisfies
\begin{itemize}
\item $a,b \in P \Rightarrow a+b \in P$
\item $a,b \in P \Rightarrow a b \in P$
\item $a \in \mathbf{R}[x_1,..,x_n] \Rightarrow a^2 \in P$
\end{itemize}
\end{definition}

\begin{definition}
Given a finite set of polynomials $S = \{a_1,...,a_m\} \subseteq  \mathbf{R}[x_1,..,x_n]$, the cone \emph{generated} by the $a_i$ is given by
\begin{equation}
\begin{aligned}
&P(S) = \\
&		\left\{p + \sum_{i=1}^r q_i b_i \ | \ p,q_1,...,q_r \in \Sigma, \ b_1,...,b_r \in M(a_i) \right\}
\end{aligned}
\end{equation}
where $\Sigma$ denotes the sum of squares cone, which is the smallest cone in $\mathbf{R}[x_1,..,x_n]$ $(\Sigma = P(\emptyset))$.
\end{definition}

Finally, we can state the Positivstellensatz, due to Stengle \cite{stengle1973nullstellensatz} (see also, e.g. \cite{bochnak1998real}). 
\begin{theorem}[\cite{stengle1973nullstellensatz}] \label{positiv}
Let $(f_i)_{i=1,...,s}$, $(g_i)_{i=1,...,t}$ and  $(h_i)_{i=1,...,u}$ be finite sets of polynomials in $\mathbf{R}[x_1,..,x_n]$. Let $P$ be the cone generated by $(f_i)_{i=1,...,s}$, $M$ the multiplicative monoid generated by $(g_i)_{i=1,...,t}$, and $I$ the ideal generated by $(h_i)_{i=1,...,u}$. Then the following are equivalent:
\begin{itemize} \label{pos}
\item $ \left\{   \begin{aligned} 
			x\in \mathbf{R}^n \left| 
				\begin{aligned} &f_i(x) \geq 0, \quad i=1,...,s  \\
						          &g_i(x) \neq 0,\quad  i=1,...,t \\ 
						          &h_i(x) = 0, \quad i=1,...,u
				\end{aligned} \right.
			\end{aligned}   \right\} = \emptyset $
\item There exist $f \in P$, $g \in M$, $h \in I$ such that $f + g^2 + h = 0$.
\end{itemize}
\end{theorem}

If the set of real solutions of the semialgebraic set is indeed empty, this result guarantees the existence of \emph{infeasibility certificates} (viz. $f \in P$, $g \in M$, $h \in I$ satisfying $f + g^2 + h = 0$); these are also known as \emph{refutations}, implying that they refute the existence of real solutions. The question then is how to search for the refutations. It turns out that the search for bounded degree refutations can be done using sum of squares and semidefinite programming, as stated by the following result from \cite{parrilo2003semidefinite}.
\begin{theorem}
Consider a semalgebraic set of the form given the Theorem \ref{pos}. Then the search for bounded degree Positivstellensatz refutations can be done using semidefinite programming. If the semialgebraic set is indeed empty and the bounded degree chosen large enough, then the semidefinite program will be feasible, and the refutations obtained from its solution. 
\end{theorem}

For a proof, see \cite{parrilo2003semidefinite}. Here, we briefly outline a constructive procedure from \cite{parrilo2003semidefinite} to find refutations. For the multiplicative monoid, let $g = \Pi_{i=1}^t g_i^{2m}$ for some $m\geq 1$. For the cone of inequalities, let $p_i$ be sum of squares polynomials of fixed degree with variable coefficients and write $$f = p_0 + p_1 f_1 + \cdots + p_s f_s  +  \cdots + p_{12...s} f_1 f_2 \cdots f_s.$$ For the ideal of equations, let $q_i$ be polynomials of fixed degree with variable coefficients and write $$ h = q_1 h_1 + \cdots q_u h_u.$$ Now consider the problem
\begin{equation}
\begin{aligned}
& \text{find} && p_i \in \Sigma, \quad q_i \in \mathbf{R}[x_1,...,x_n] \\
& \text{subject to} && f+g^2+h=0 \\
\end{aligned}
\end{equation}
where the decision variables are the coefficients of the $p_i$ and $q_i$. For a given degree of the $p_i$ and $q_i$ this is a semidefinite programming feasibility problem, with the Positivstellensatz constraint leading to affine constraints on the coefficients of the $p_i$ and $q_i$. If the semialgebraic set is indeed empty, then there exist polynomial Positivstellensatz certificates according to Theorem \ref{positiv}. By construction of the semidefinite program above, there is a finite number $d^*$ such that for every $d \geq d^*$ the program is feasible, and the certificates can be obtained directly from a feasible solution of the program. 

For computational complexity reasons, there is no guarantee of obtaining low degree certificates for every possible instance. However, low degree certificates have been obtained for many practical problems \cite{parrilo2003semidefinite}.

\section{Sum of Squares Programming and Distance-Based Formation Shape Control}
We now apply sum of squares programming and the Positivstellensatz to a distance-based formation shape control problem. We will show that the question of existence of a locally stable incorrect equilibrium shape corresponds to the question of whether or not a certain semialgebraic set, consisting of finite sets of polynomial equations, inequations, and inequalities, is empty. By the arguments in the previous section, this can be checked in principle using semidefinite programming. In particular, for a given desired formation shape, if there does not exist a locally stable incorrect equilibrium shape, one can obtain a \emph{certificate} via the Positivstellensatz.

\subsection{Equations of Motion}
To simplify the exposition, we will consider the equations of motion for a four-agent formation shape control problem in the plane; however, the method extends to any size formation and to higher dimensions in principle. 

Let $p=[p_1^T,p_2^T,p_3^T,p_4^T]^T\in \mathbf{R}^8$ be a vector of the four agent
positions in the plane. We use a single
integrator agent dynamic model for each agent
\begin{equation}
\dot{p}_i=u_i \nonumber
\end{equation}
where $u_i$ is the control input to be specified. Let
$\bar{d}=[\bar{d}_{12},\bar{d}_{13},\bar{d}_{14},\bar{d}_{23},\bar{d}_{24},\bar{d}_{34}]^T$ be a vector of
desired interagent distances that define the formation shape and are
to be actively controlled. We assume that the entries of $\bar{d}$
correspond to a realizable shape. Let $d(p)=[d_{12}(p),d_{13}(p),d_{14}(p),d_{23}(p),d_{24}(p),d_{34}(p)]^T$ denote instantaneous interagent distances, i.e.
\begin{equation} \nonumber
\begin{split}
d^2(p)=[&||p_1-p_2||^2,||p_1-p_3||^2,||p_1-p_4||^2,\\&||p_2-p_3||^2,||p_2-p_4||^2,||p_3-p_4||^2]^T
\end{split}
\end{equation}
The distance error function is given by
\begin{equation} \label{errorfunction}
\begin{split}
e(p)&=d^2(p)-\bar{d}^2\\&=[e_{12},e_{13},e_{14},e_{23},e_{24},e_{34}]^T
\end{split}
\end{equation}

Now consider the potential function
\begin{equation} \nonumber
\begin{split}
V(p)&=\frac{1}{2}||e(p)||^2\\&=\frac{1}{2}(e_{12}^2+e_{13}^2+e_{14}^2+e_{23}^2+e_{24}^2+e_{34}^2),
\end{split}
\end{equation} which quantifies the total interagent
distance error between the current formation and the desired
formation $\bar{d}$. Following Krick et al. \cite{krick2009stabilisation}, we define the gradient control law
\begin{equation} \label{gradientcontrol}
u=-\nabla V(p)^T,
\end{equation}
which yields the closed-loop system
\begin{eqnarray} \nonumber \label{closedloop}
\dot{p}&=&-\nabla V(p)^T \\
       &=&-[J_r(p)]^Te(p) \\
       &=&-(E(p)\otimes I_2) p \nonumber
\end{eqnarray}
where $J_r(p)$ is the Jacobian of the distance error function $e(p)$ (this matrix is also called the \emph{rigidity matrix} in rigid graph theory), $\otimes$ is the Kronecker product and the matrix $E(p)$ is given by 
\begin{equation} \label{ematrix} 
\begin{aligned}
&E(p)= \left[
        \begin{array}{cccc}
          \sum_{j} e_{1j} & -e_{12} & -e_{13} & -e_{14} \\
          -e_{12} & \sum_j e_{2j} & -e_{23} & -e_{24} \\
          -e_{13} & -e_{23} & \sum_j e_{3j} & -e_{34} \\
          -e_{14} & -e_{24} & -e_{34} & \sum_j e_{4j} \\
        \end{array}
      \right]
\end{aligned}
\end{equation} 
with $e_{ij}=e_{ji}$ and $e_{ii} = 0$. Agent-wise, the equations of motion are
\begin{equation} \label{agentequil}
\begin{split}
\dot{p}_1 &= e_{12}(p_2-p_1) + e_{13}(p_3-p_1) + e_{14}(p_4-p_1) \\
\dot{p}_2 &= e_{12}(p_1-p_2) + e_{23}(p_3-p_2) + e_{24}(p_4-p_2) \\
\dot{p}_3 &= e_{13}(p_1-p_3) + e_{23}(p_2-p_3) + e_{34}(p_4-p_3) \\
\dot{p}_4 &= e_{14}(p_1-p_4) + e_{24}(p_2-p_4) + e_{34}(p_3-p_4)
\end{split}
\end{equation}

The equilibrium points\footnote{The equilibrium points form sets, with each set obtainable by translation and rotation of any point in the set.} of
the closed-loop system (\ref{closedloop}) are the
critical points of the potential function $V$. The Jacobian of the
right side of (\ref{closedloop}) is
the same as the negative of the Hessian of $V$, which is given by
\begin{equation} \label{Hessian}
H_V(p)= 2 J_r(p)^T J_r(p) + E(p)\otimes I_2
\end{equation}
The stability of any one equilibrium point and its associated set of equilibrium points of
(\ref{closedloop}) are determined by the nature of the critical
points of $V$: local minima are locally stable and local maxima and saddle
points are locally unstable (noting that for all equilibria, there will always be three eigenvalues of $H_V$ which are zero, corresponding to the rotational and translational invariance). 

\subsection{Applying the Positivstellensatz}
To apply the Positivstellensatz, observe that the conditions for equilibrium, incorrectness and local stability are all given by finite sets of polynomials in $p$. Equilibrium is given by the polynomial equations $(E(p)\otimes I_2) p = 0$, a cubic function of $p$. So let $h_i(p)$ ($i=1.,,,.8$) be the rows of $(E(p)\otimes I_2) p$. Shape incorrectness is equivalent to the polynomial inequation of the error vector $\Vert e(p) \Vert^2 \neq 0$, so let $g_i(p)=\Vert e(p) \Vert^2$ ($i = 1$). Local stability requires positive semidefiniteness of $H_V(p)$, which is equivalent to the non-negativity of all principal minors of $H_V(p)$, so let $f_i(p)$ ($i=1,...,2^8-1$) be the principal minors of $H_V(p)$. For a given desired formation $\bar{d}$, the question of existence of locally stable incorrect equilibria is the same as asking if the semialgebraic set
 $$ \left\{   \begin{aligned} 
			p\in \mathbf{R}^8 \left| 
				\begin{aligned} &f_i(p) \geq 0, \quad i=1,...,2^8-1  \\
						          &g_i(p) \neq 0,\quad  i=1,...,6 \\ 
						          &h_i(p) = 0, \quad i=1,...,8
				\end{aligned} \right.
			\end{aligned}   \right\}  $$
is empty. Thus, one can use the semidefinite programming techniques described in the previous section to certify non-existence of locally stable incorrect equilibria. 

For larger formations (including in higher dimensions and with any graph underlying the information architecture), the procedure to construct the above semialgebraic set is the same. The conditions for equilibrium, shape incorrectness, and local stability can all still be written as finite sets of polynomials, with the only difference being that the number of variables and the number of polynomials increases. Thus, this approach extends to any size formation and higher dimensions in principle.
\begin{remark} \nonumber
Again, if the set is indeed empty, only \emph{existence} of certificates is guaranteed; one cannot always expect to find low degree certificates. The size of the semidefinite program grows very fast with increasing number of variables and degree of the polynomial certificates, so the method may quickly become computationally prohibitive for larger formations. Since the conditions are independent of rotation and translation, one can fix e.g. $p_1=0$ and $p_{2x}=0$, thereby reducing the number of variables by three. 
\end{remark}

\subsection{Simplifications and Extensions}
It is possible to reduce the computational complexity of the search for certificates by using simpler (but more conservative) necessary conditions for equilibrium, shape incorrectness, or local stability. There is perhaps most to gain from replacing the local stability condition with a simpler necessary condition. The principal minor test for positive semidefiniteness involves checking all principal minors, the number of which is exponential in the number of agents, so the number of $f_i$'s grows very fast. We might hope still to obtain certificates by replacing the principal minor test with a more conservative necessary condition for positive semidefiniteness, e.g. only check some principal minors. Specifically, it is possible to consider positive definiteness of a suitably selected $5 \times 5$ principal submatrix of $H_V(p)$. Fixing $p_1=0$ and $p_{2x}=0$, then the last 5 by 5 rows and columns of $H_V(p)$ must be positive definite, which only involves checking $5$ principal minors.

For the equilibrium and shape incorrectness conditions, one could consider a collections of smaller alternative problems involving only the individual entries of $(E(p)\otimes I_2) p = 0 $ and individual entries of $e(p) \neq 0$, and combinations thereof.

A further difficulty with the local stability condition is that the principal minors involve very large expressions for determinants of submatrices of $H_V(p)$. One requires a computer to form and parse these expressions to obtain constraints for the semidefinite program.

The method as stated can determine non-existence of a locally stable incorrect equilibrium shape for a \emph{given} desired formation shape. It is also straightforward in principle to consider the desired formation $\bar{d}$ as a variable: the $f_i$'s, $g_i$'s and $h_i$'s are now all also polynomials in $\bar{d}$. One would then have an additional equality constraint $k(\bar{d})$ corresponding to, for example, $\bar{d}$ having Cayley-Menger determinant equal to zero \cite{blumenthal1970theory}:
\begin{equation}
k(\bar{d}) = \left|\begin{array}{ccccc}0 & \bar{d}_{12}^{2} & \bar{d}_{13}^{2} & \bar{d}_{14}^{2} & 1 \\ \bar{d}_{12}^{2} & 0 & \bar{d}_{23}^{2} & \bar{d}_{24}^{2} & 1 \\ \bar{d}_{13}^{2} & \bar{d}_{23}^{2} & 0 & \bar{d}_{34}^{2}  & 1 \\ \bar{d}_{14}^{2} & \bar{d}_{24}^{2} & \bar{d}_{34}^{2} & 0 & 1 \\1 & 1 & 1 & 1 & 0\end{array}\right| = 0.
\end{equation}
It is possible to work with other similar and possible similar constraints involving Gramian-like matrices; see \cite{dasgupta2011controlling}. 
Then the question of non-existence of a locally stable incorrect formation for \emph{any} desired four-agent formation is the same as asking if the semialgebraic set
 $$ \left\{   \begin{aligned} 
			p\in \mathbf{R}^8, \bar{d}\in \mathbf{R}^6 \left| 
				\begin{aligned} &f_i(p,\bar{d}) \geq 0, \quad i=1,...,2^8-1  \\
						          &g_i(p,\bar{d}) \neq 0,\quad  i=1,...,6 \\ 
						          &h_i(p,\bar{d}) = 0, \quad i=1,...,8 \\
						          &k(\bar{d}) = 0
				\end{aligned} \right.
			\end{aligned}   \right\}  $$
is empty.

\section{Conclusions}
This paper has shown that sum of squares programming techniques and the Positivstellensatz provides a method in principle to certify non-existence of locally stable incorrect equilibrium shapes in a distance-based formation shape control problem. The certificates can be obtained by solving an explicitly defined semidefinite program, for which well-developed software packages exist. The method extends to any size formation and in higher dimension in principle, though the computational complexity grows very fast with the problem size.

Preliminary computational results have not yet yielded certificates, though this is an immediate direction for future work. For certain formation shape control problems, there may very well exist locally stable incorrect equilibrium shapes. These problems may require alternative approaches to guarantee global stability, such as using a combination of distances and angles. Preliminary work in this direction can be found in \cite{bishop2012control} and \cite{bishop2013control}. The ultimate goal is to develop methods to control an arbitrary formation to a uniquely specified shape.

\bibliographystyle{plain}  
\bibliography{refs}  

\end{document}